\documentclass[12pt]{article}

\usepackage{amssymb}
\usepackage{amsmath}

 \newtheorem{thm}{Theorem}[section]
 \newtheorem{prop}{Proposition}[section]
 \newtheorem{lem}{Lemma}[section]
 \newtheorem{cor}{Corollary}[section]
 \newtheorem{exm}{Example}[section]
 \newtheorem{dfn}{Definition}[section]
 \newtheorem{rem}{Remark}[section]

\renewcommand{\leq}{\leqslant}
\renewcommand{\geq}{\geqslant}

\newcommand{\rbx}{\hfill{\rule{1ex}{1ex}}}

\newcommand{\ind}{\mathrm{ind}\,}
\newcommand{\Span}{\mathrm{span}\,}
\newcommand{\diag}{\mathrm{diag}\,}

 \newcommand{\vp}{\varphi}

\newcommand{\coker}{\mbox{\rm coker}\,}

\newcommand{\im}{\mathrm{im}\,}
\newcommand{\sign}{\mathrm{sgn}\,}

\newcommand{\nn}{\nonumber}

 \newcommand{\esssup}{\hbox{\rm ess}\sup_{\!\!\!\!\!\!\!\!\!\!\! t\in {\mathbb{T}}}}

\newcommand{\cJ}{\mathcal{J}}

\newcommand{\cP}{\mathcal{P}}
\newcommand{\cQ}{\mathcal{Q}}
\newcommand{\cR}{\mathcal{R}}

\newcommand{\sC}{{\mathbb C}}

\newcommand{\sN}{{\mathbb N}}

\newcommand{\sT}{{\mathbb T}}

\newcommand{\sZ}{{\mathbb Z}}

\begin{document}

\vspace*{10mm}

\begin{center}
{\Large\textbf{Closed form solution of non-homogeneous\\[1ex] equations
with Toeplitz plus Hankel operators}}
\end{center}

\vspace{5mm}

\begin{center}

\textbf{Victor D. Didenko and Bernd Silbermann}

\vspace{2mm}


Universiti Brunei Darussalam, Bandar Seri Begawan, BE1410  Brunei;
diviol@gmail.com

Technische Universit{\"a}t Chemnitz, Fakult{\"a}t f\"ur Mathematik,
09107 Chemnitz, Germany; silbermn@mathematik.tu-chemnitz.de

 \end{center}


  \vspace{10mm}



\begin{abstract}
Considered is the equation
 $$
 (T(a)+H(b))\phi=f,\eqno(\star)
 $$
where $T(a)$ and $H(b)$, $a,b\in L^\infty(\mathbb{T})$ are,
respectively, Toeplitz and Hankel operators  acting on the classical
Hardy spaces $H^p(\mathbb{T})$, $1<p<\infty$. If the generating
functions $a$ and $b$ satisfy the so-called matching condition
[1,2],
 $$
a(t) a(1/t)=b(t)b(1/t), \, t\in \mathbb{T},
 $$
an efficient method for solving equation ($\star$) is proposed. The
method is based on the Wiener--Hopf factorization of the scalar
functions $c(t)=a(t)b^{-1}(t)$ and $d(t)=a(t)b^{-1}(1/t)$ and allows
one to find all solutions of the equations mentioned.
 \end{abstract}

\section{Introduction\label{s1}}
Let $\sT:=\{ t\in\sC:|t|=1 \}$ be the counterclockwise oriented unit
circle in the complex plane $\sC$, and let $L^p=L^p(\sT)$, $1\leq
p\leq \infty$ denote the space of all Lebesgue measurable functions
$f$ such that
 \begin{align*}
||f||_p &: = \left ( \int_\sT |f(t)|^p \,| dt| \right )^{1/p}
<\infty,\quad  1\leq p <\infty,\\
 ||f||_\infty&: = \esssup |f(t)|<\infty.
\end{align*}
If $f\in L^1$ then, as usual, we use the notation $\widehat{f}_n$ to
denote the $n$-th Fourier coefficient of $f$, i.e.
 \begin{equation*}
\widehat{f}_n:=\frac{1}{2\pi}\int_0^{2\pi}f(e^{i\theta})
e^{-in\theta}\, d\theta\, , \quad n\in\sZ.
 \end{equation*}
The Hardy spaces $H^p:=H^p(\sT)$ and
$\overline{H^p}:=\overline{H^p(\sT)}$, $1\leq p \leq \infty$ are
defined as follows
\begin{equation*}
 \begin{aligned}
   H^p:=\{ f\in L^p: \widehat{f}_n=0 \,\text{ for all } n <0\}, \\
 \overline{H^p}:=\{ f\in L^p: \widehat{f}_n=0 \,\text{ for all } n
 >0\}.
\end{aligned}
\end{equation*}
Let $\sZ_+:=\sN\cup \{0\}$ denote the set of all non-negative
integers. Consider the operator $P$ defined by
 \begin{equation*}
 P:  \sum_{n\in \sZ} \widehat{f}_n t^n \to \sum_{n\in \sZ_+}^\infty
 \widehat{f}_n t^n \, ,\\
 \end{equation*}
and let $Q:=I-P$, where $I$ is the identity operator. It is easily
seen that $P$ and $Q$ are complimentary projections and it is
well-known that they are bounded on any space $L^p$,
$p\in(1,\infty)$. Along with the operators $P$ and $Q$ we also
consider the operator $J:L^p\to L^p$,
\begin{equation*}
 (Jf)(t) := t^{-1}f(t^{-1})\,, \quad t\in\sT .
 \end{equation*}

For $a,b\in L^\infty$, the Toeplitz plus Hankel operator
$T(a)+H(b):H^p\to H^p$ is defined by
  \begin{equation}\label{eqn1}
    T(a)+H(b):= PaP + PbQJ.
\end{equation}
Various properties of these operators have been established in
literature. Special interest represents the Fredholmness of the
operator $T(a)+H(b)$ and investigation of its kernel and cokernel.
In the case of piecewise continuous generating functions $a$ and
$b$, Fredholm properties of the operator \eqref{eqn1} can be derived
by a direct application of results \cite[Sections
4.95-4.102]{BS:2006}, \cite[Sections 4.5 and 5.7]{RSS:2011},
\cite{RS1990}. The case of quasi piecewise continuous generating
functions has been studied in \cite{Si:1987}, whereas formulas for
the index of the operators \eqref{eqn1}, considered on different
Banach and Hilbert spaces and with various assumptions about the
generating functions $a$ and $b$, have been established in
\cite{DS:2013,RS:2012}. Recently, progress has been made in
computation of defect numbers $\dim\ker (T(a)+H(b))$ and $\dim\coker
(T(a)+H(b))$ for various classes of generating functions $a$ and $b$
\cite{BE:2013,DS:2014a}. The more delicate problem of the
description of the spaces $\ker (T(a)+H(b))$ and $\coker
(T(a)+H(b))$ has been considered \cite{DS:2014a, DS:2014}. In
particular, for generating functions $a$ and $b$, satisfying the
so-called matching condition (see condition \eqref{eqn12} below),
explicit and efficient formulas for the elements of the kernel of
the operator $T(a)+H(b)$ have been established. Thereby, for this
class of generating functions, homogeneous equations with Toeplitz
plus Hankel operators can be effectively solved. On the other hand,
non-homogeneous equations with Toeplitz plus Hankel operators and
also equations with Wiener--Hopf plus Hankel operators often arise
in applications \cite{Du:1979, Karapetiants2001, MST:1992} but, as a
rule, only some approximation methods of their solutions have been
systematically studied so far.

The aim of this work is to present a method for solution of the
operator equations
  \begin{equation}\label{eqn2}
  (T(a)+H(b))\vp =f, \quad f\in L^p, \quad 1<p<\infty,
 \end{equation}
with generating functions $a,b\in L^\infty$ satisfying relation
\eqref{eqn12}. Note that this relation has been first used in
\cite{BE:2013} when studying the dimension of kernels and cokernels
of Toeplitz plus Hankel operators with piecewise continuous
generating functions $a$ and $b$. Regardless of \cite{BE:2013}, the
importance of relation \eqref{eqn12} for the investigation of
Toeplitz plus Hankel operators has been mentioned in \cite[Remark
9]{DS:2013}.

The approach proposed is also applicable to non-homogeneous
equations with Wiener--Hopf plus Hankel operators with generating
matching functions. However, at the moment there is no description
for the kernels of such operators. Nevertheless, an investigation of
Wiener--Hopf plus Hankel operators in the situation mentioned, has
been started in \cite{DS:2014b}. As soon as that work will be
completed, relevant results can be also used for solution of
non-homogeneous Wiener--Hopf plus Hankel equations. The details of
the corresponding study will be presented elsewhere.

\section{Toeplitz plus Hankel equations and equations with matrix Toeplitz
operators \label{s2}}

In this section we establish connections between the solutions of
Toeplitz plus Hankel equation \eqref{eqn2} and solutions of an
equations with a matrix Toeplitz operator. On the space $H^p\times
H^p$, $1<p<\infty$, let us consider the operator $\cR$,
  \begin{equation*}
   \cR:= \diag(T(a)+H(b),T(a)-H(b))\, ,
 \end{equation*}
and let $\cP$ and $\cQ$ denote the operators
 $$
\cP:=\diag(P,P), \quad \cQ:=\diag(Q,Q)\, ,
  $$
acting on the space $L^p(\sT)\times L^p(\sT)$, $1<p<\infty$. For any
element $a\in L^\infty$, set $\widetilde{a}(t):=a(1/t)$. Let
$C=C(a,b)$ and  $V=V(a,b)$ be $2\times 2$-matrices,
 \begin{align}
 & C(a,b) :=\left(%
\begin{array}{cc}
  1 & 0 \\
 \widetilde{b} & \widetilde{a} \\
   \end{array}%
\right),
 \quad
  V(a,b) :=\left(%
\begin{array}{cc}
  a-b \widetilde{b} \widetilde{a}^{-1} & b \widetilde{a}^{-1} \\
 -\widetilde{b}\widetilde{a}^{-1} & \widetilde{a}^{-1} \\
   \end{array}%
\right), \nn\\
\intertext{and let $\cJ,A_1,  A_2, B, R$ be the operators defined by
}
 & \cJ:=  \frac{1}{2}\left(%
\begin{array}{cc}
  I & J \\
 I & -J \\
\end{array}
\right) , \quad A_1:=\diag(I,I)-\diag(P,Q) \left(%
\begin{array}{cc}
  a & b \\
  \widetilde{b}  & \widetilde{a}  \\
   \end{array}%
\right) \diag(Q,P), \nn\\[1ex]   
& A_2:=\diag(I,I)+\cP V(a,b) \cQ, \nn\\[1ex]
&  B:=\cP\,V(a,b)\,\cP+\cQ, \nn\\[1ex]
& R:= \cR+\cQ. \nn
  \end{align}
According to the relation (3.4) of \cite{DS:2014}, the operators $R$
and $B$ are connected as follows
   \begin{equation}\label{eqn3}
    \cJ^{-1}  R  \cJ= A_1  A_2  B  C.
 \end{equation}

  \begin{lem}\label{l1}
Let $y=(g,h)^T \in H^p(\sT) \times H^p(\sT)$. If\/
$x_0=(\vp,\psi)^T\in H^p(\sT) \times H^p(\sT)$ is a solution of the
equation
   \begin{equation}\label{eqn4}
   \cR x =y,
 \end{equation}
then the element
  \begin{equation}\label{eqn5}
X_0= \left (
\vp+\psi,P(\widetilde{b}(\vp+\psi)+\widetilde{a}J(\vp-\psi))
\right)^T
 \end{equation}
is a solution of the equation
  \begin{equation}\label{eqn6}
  \cP V(a,b)\cP X= Y,
 \end{equation}
 where $Y= \cP A_2^{-1}A_1^{-1}\cJ^{-1}y$.
   \end{lem}
\textbf{Proof.} If $x_0=(\vp,\psi)^T$ is a solution of \eqref{eqn4},
then one also has
  $$
(\cR + \cQ) x_0=y.
  $$
Consequently,
  $$
\cJ^{-1}(\cR + \cQ)\cJ \cJ^{-1}x_0=\cJ^{-1} y,
  $$
and relation \eqref{eqn3} implies that
 $$
BC\cJ^{-1} x_0= A_2^{-1} A_1^{-1} \cJ^{-1} y,
 $$
or
 \begin{equation*}
  (\cP\,V(a,b)\,\cP+\cQ)  C \cJ^{-1}
x_0=   A_2^{-1} A_1^{-1}\cJ^{-1}y.
 \end{equation*}
Therefore, the element $X_0=\cP C \cJ^{-1} x_0$ is a solution of the
equation \eqref{eqn6} and since
 $$
\cJ^{-1}= \left(%
\begin{array}{cc}
  I & I \\
 J & -J \\
\end{array}
\right),
 $$
the relation \eqref{eqn5} follows.
 \rbx

 Thus the solutions of the equation \eqref{eqn4} generate solutions
of equation \eqref{eqn6} with the corresponding right-hand side $Y$.
On the other hand, a similar reverse statement is also true.
 \begin{lem}\label{l2}
If the equation \eqref{eqn6} is solvable and if $X_0=(\Phi,\Psi)$ is
a solution of \eqref{eqn6} with the right-hand side $Y=(g,h)\in
H^p\times H^p$, then the element
 \begin{equation}\label{eqn7}
 x_0=\frac{1}{2}(\Phi-JQc\Phi+JQ\widetilde{a}^{-1}\Psi,
      \Phi+JQc\Phi-JQ\widetilde{a}^{-1}\Psi)^T,
 \end{equation}
where $c=\widetilde{b}\widetilde{a}^{-1}$, is a solution of the
equation \eqref{eqn4} with the right-hand side
   \begin{equation}\label{eqn8}
 y= \cP\cJ A_1 A_2 Y.
 \end{equation}
  \end{lem}

   \textbf{Proof.}
If $X_0$ is a solution of equation \eqref{eqn6}, then
 $$
( \cP V(a,b)\cP +\cQ)X_0= Y.
 $$
Applying the operator $\cJ A_1 A_2$ to this equation, one obtains
 $$
\cJ A_1 A_2( \cP V(a,b)\cP +\cQ) X_0= \cJ A_1 A_2 Y.
 $$
The last expression can be rewritten as
  $$
(\cJ A_1 A_2( \cP V(a,b)\cP +\cQ)C \cJ^{-1}) \cJ C^{-1} X_0= \cJ A_1
A_2 Y,
 $$
and the relation \eqref{eqn3} implies that
 $$
R \cJ C^{-1} X_0= \cJ A_1 A_2 Y.
 $$
Therefore, the element $z_0$,
 $$
z_0= \cJ C^{-1} X_0=\frac{1}{2}(\Phi-Jc\Phi+J\widetilde{a}^{-1}\Psi,
      \Phi+Jc\Phi-J\widetilde{a}^{-1}\Psi)^T
 $$
satisfies the equation
 $$
Rx=\cJ A_1 A_2 Y,
 $$
and since
 $$
Rz_0=(\cR+\cQ) z_0 = \cJ A_1 A_2 Y= \cP\cJ A_1 A_2 Y + \cQ\cJ A_1
A_2 Y,
 $$
the element $x_0=\cP z_0$ is a solution of equation \eqref{eqn4}.
\rbx

As we will see later on, in some cases the equation \eqref{eqn6} can
be resolved in a closed form. Thus using the representation
\eqref{eqn7} one can also derive solutions $\vp$ of \eqref{eqn2}
from the solutions $(\Phi, \Psi)$ of the equation \eqref{eqn6},
provided that the right-hand side $Y\in H^p\times H^p$ in
\eqref{eqn6} is chosen in such a way that the first coordinate of
the vector $y$ in \eqref{eqn8} is equal to the right hand side $f$
of the equation \eqref{eqn2}. In this case, a solution of the
equation \eqref{eqn2} is given by the formula
  \begin{equation}\label{eqn9}
 \vp_0=\frac 12 (\Phi-JQc\Phi+JQ\widetilde{a}^{-1}\Psi).
 \end{equation}
It is clear that the choice of possible right-hand side $Y$ in
\eqref{eqn6} is not unique. Nevertheless, it seems that the most
suitable initial vector $Y$ in the equation \eqref{eqn6} is the one
transferred into the vector $(f,0)^T$ by the operator $\cP \cJ A_1
A_2$. Indeed, with such a choice of $Y$, the second equation in
\eqref{eqn4} would have the form
 $$
(T(a)-H(b)) \psi=0.
 $$
This equation is always solvable, so there will be no additional
condition related to the solvability of the equation \eqref{eqn4}.
In order to find such a right-hand side $Y$, one has to resolve the
operator equation
 \begin{equation}\label{eqn10}
 \cP \cJ A_1 A_2 Y=(f,0)^T.
 \end{equation}
Setting $Y=(g,h)^T$, $g,h\in H^p$ one can write equation
\eqref{eqn10} as the system of equations
 \begin{equation}\label{eqn11}
 \left \{
  \begin{aligned}
    g-PbPh-PaQJh&=2f\\
    g-PbPh+PaQJh&=0,
\end{aligned}
\right.
 \end{equation}
with respect to unknown functions $g$ and $h$. However, the problem
of the determination of the functions $g$ and $h$ from the system
\eqref{eqn11} is equivalent to the solution of the equation
 $$
PaQJh=-2f,
 $$
which is not a simple task. Nevertheless, system \eqref{eqn11}
suggests a simple choice of the functions $g$ and $h$ which would
lead to the equation \eqref{eqn2} with the required right-hand side
$f$, viz. one can consider the pair $g=2f$, $h=0$ so the
corresponding right-hand side $Y=(2f,0)$. Of course, a consequence
of such a choice of the right hand side $Y$ is that the equation
 $$
(T(a)-H(b))\psi=f
 $$
must be also solvable. Note that the solvability of the last
equation is not directly connected to the solvability of
\eqref{eqn2}, which means that if the operator $\cP V(a,b)\cP$ is
not right invertible, our method will not work for some right hand
sides $f$ from the image of the operator $T(a)+H(b)$. However, the
set of acceptable right-hand sides $f$ is quite large because it is
generates by those pairs $(f,0)^T$, $f\in H^p$ which belongs to the
image of the operator $\cP V(a,b) \cP$. More precisely, the
following proposition holds.

 \begin{prop}\label{p3}
If $Y:=(2f,0)^T$, $f\in H^p$ and the equation \eqref{eqn6} with the
right-hand side $Y$ is solvable with the solution $(\Phi,\Psi)$,
then the equation \eqref{eqn2} is also solvable and one of its
solution can be written in the form \eqref{eqn9}.
 \end{prop}
 \textbf{Proof.}
 Straightforward computation.
\rbx

Note that the connection between the solvability of the systems
\eqref{eqn4} and \eqref{eqn6} with the corresponding right-hand
sides $(f,f)^T$ and $(2f,0)^T$ will be discussed later (see Remark
\ref{rem1} below).

\section{Solution of non-homogeneous equations with\\ Toeplitz plus Hankel operators.}

In this section we construct solutions of the non-homogeneous
equation \eqref{eqn6} in the case where the generating functions
$a,b \in L^\infty$ are connected in a special way. Thus let us
assume that $a$ and $b$ satisfy the relation
 \begin{equation}\label{eqn12}
 a(t)\widetilde{ a}(t)=b(t) \widetilde{b}(t), \quad t\in \sT,
 \end{equation}
where $\widetilde{a}(t)=a(1/t)$ and $\widetilde{b}(t)=b(1/t)$, as
before. In what follows this relation is called the matching
condition and any duo $(a,b)$ with the property \eqref{eqn12} is
called the matching pair. For any matching pair $(a,b)$ one can
construct another pair $(c,d)$ of functions $c$ and $d$ defined by
  $$
c(t):=ab^{-1}(=\widetilde{b}\widetilde{a}^{-1}), \quad d:= b
\widetilde{a}^{-1}(=\widetilde{b}^{-1} a).
  $$
Such a pair $(c,d)$ is called the subordinated pair for $(a,b)$, and
the functions which constitutes a subordinated pair possess the
property
 $$
c\widetilde{c}=1=d\widetilde{d}.
 $$
Let us also point out that if $(c,d)$ is the  subordinated  pair for
a matching pair $(a,b)$, then $(\overline{d},\overline{c})$ is the
subordinated pair for the matching pair $(\overline{a},
\widetilde{\overline{b}})$ which define the adjoint operator
 \begin{equation*}
(T(a)+H(b))^* :=T(\overline{a})+H(\widetilde{\overline{b}}),
 \end{equation*}
for the Toeplitz plus Hankel operator $T(a)+H(b)$.

Consider now the equation
  \begin{equation}\label{eqn13}
  T(V(a,b))X=(2f,0)^T, \quad f\in H^p,
 \end{equation}
where $T(V(a,b))=\cP V(a,b) \cP$ is the Toeplitz operator generated
by the matrix $V(a,b)$. If $(a,b)$ is a matching pair, then $V(a,b)$
is a triangular matrix, namely,
 $$
V(a,b)=\left(%
\begin{array}{cc}
 0  & d \\
  -c  & \widetilde{a}^{-1} \\
   \end{array}%
\right),
 $$
where $(c,d)$ is the corresponding subordinated pair.

It turns out that in this case, the solutions of the matrix equation
\eqref{eqn13} depend on the solutions of equations with Toeplitz
operators $T(c)$ and $T(d)$. Therefore, we would like to recall some
facts concerning scalar Toeplitz operators.
 \begin{dfn}
A function $a\in L^\infty$ admits a week Wiener--Hopf factorization
in $H^p$, if it can be represented in the form
 \begin{equation}\label{eqn14}
    a=a_- t^n a_+,
\end{equation}
where $n\in \sZ$, $a_+\in H^q$, $a_+^{-1}\in H^p$, $a_-\in
\overline{H^p}$, $a_-^{-1}\in \overline{H^q}$, and $a_-(\infty)=1$.
 \end{dfn}
The weak Wiener--Hopf factorization of a function $a$ is unique, if
it exists.  The functions $a_-$ and $a_+$ are called the
factorization factors, and the number $n$ is the factorization
index. If $a\in L^\infty$ and the operator $T(a)$ is Fredholm, then
the function $a$ admits the weak Wiener--Hopf factorization with
$n=-\ind T(a)$ \cite{BS:2006, LS:1987}. Moreover, in this case, the
factorization factors possess an additional property--viz. the
linear operator $a_+^{-1}Pa_-^{-1}I$ defined on $\Span
\{t^k:k\in\sZ_+\}$ can be boundedly extended on the whole space
$H^p$. Throughout this paper, such a kind of the weak Wiener--Hopf
factorization in $H^p$ is called simply Wiener--Hopf factorization
in $H^p$. The following result is well--known.
 \begin{thm}[see {\cite{BS:2006}}]\label{t11}
If $a\in L^\infty$, then Toeplitz operator $T(a):H^p\to H^p$,
$1<p<\infty$ is Fredholm and $\ind T(a)=-n$ if and only if the
generating function $a$ admits the Wiener--Hopf factorization
\eqref{eqn14} in $H^p$.
 \end{thm}
Thus if $T(a)$ is Fredholm and $n\geq 0$, then the operator $T(a)$
is left invertible and
  \begin{equation*}
T_l^{-1}(a)= T(t^{-n}) T(a_+^{-1}) T(a_-^{-1})
 \end{equation*}
is a left inverse of $T(a)$, whereas for $n\leq 0$, the operator
$T(a)$ is right invertible and
  \begin{equation}\label{eqn15}
T_r^{-1}(a)= T(a_+^{-1}) T(a_-^{-1})T(t^{-n})
 \end{equation}
is one of its right inverses. For definiteness, in the following the
notation $T_r^{-1}(a)$ is always used for the operator defined by
\eqref{eqn15}.

In this work we mainly use Toeplitz operators generated by matching
functions, i.e. by the functions satisfying the condition $g
\widetilde{g}=1$. For such functions, factorization \eqref{eqn14}
has special properties. Thus if $T(g):H^p\to H^p$, $1<p<\infty$ is a
Fredholm operator, then it was shown in \cite{DS:2014} that the
matching function $g$ can be represented as follows
  \begin{equation}\label{Eq27}
  g(t)= g_+(t) t^n \left (\boldsymbol\sigma(g) \widetilde{g}_+^{-1}(t) \right ),
 \end{equation}
where $\boldsymbol\sigma(g):=g_+(0)=\pm 1$ is called the
factorization signature. In \cite{DS:2014}, there is discussed how
to find the factorization signature in some particular cases. For
instance, if $g$ is continuous at the point $1$ then
$\boldsymbol\sigma(g)=g(1)$.

Assume now that the operators $T(c)$ and $T(d)$ are Fredholm and let
$\kappa_1=\kappa_c:=\ind T(c)$ and $\kappa_2=\kappa_d:=\ind T(d)$.
If $T(c)$ is a right invertible operator, then we also consider the
operator $W$ defined by
  \begin{equation*}
  W\vp := T_r^{-1}(c)T(\widetilde{a}^{-1}) \vp - JQcP
   T_r^{-1}(c)T(\widetilde{a}^{-1})\vp + JQ \widetilde{a}^{-1} \vp.
 \end{equation*}
Now we can describe the solution set for the equation \eqref{eqn2}.
The structure of this set depends on the indices of the operators
$T(c)$ and $T(d)$. Let us start with the case where $\kappa_c$ and
$\kappa_d$ are non-negative.
 \begin{thm}\label{t1}
If $\kappa_c\geq 0$ and $\kappa_d\geq0$, then for any $f\in H^p$ the
equation \eqref{eqn2} is solvable and all solutions of this equation
are given by the formula
 \begin{multline}\label{eqn16}
\vp =T_r^{-1}(c) T(\widetilde{a}^{-1})T_r^{-1}(d) f -
JQcT_r^{-1}(c)T(\widetilde{a}^{-1})T_r^{-1}(d) f +JQ \widetilde{a}^{-1} T_r^{-1}(d) f \\
  \quad + \sign(\kappa_c)\, c_+^{-1}\sum_{k=0}^{\kappa(1)} r_k^{(1)} u_k^{(\kappa(1),-)} +
   \sign(\kappa_d)\, W \left (d_+^{-1}\sum_{k=0}^{\kappa(2)} r_k^{(2)} u_k^{(\kappa(2),+)}\right ),
\end{multline}
where $r_k^{(i)}$, $k=0,1, \ldots \, , k(i)$, $i=1,2$, are arbitrary
complex numbers,
 \begin{equation*}
  \begin{aligned}
&\kappa(i)=m_i-1, \; u_k^{(\kappa(i),\pm)}(t) = t^{m_i-k-1} \pm
\boldsymbol\sigma(s_i)t^{m_i+k},  \text{if} \quad \kappa_i=2m_i, \;
i=1,2, \\
 &\kappa(i)=m_i, \quad u_k^{(\kappa(i),\pm)}(t) = t^{m_i+k} \pm
\boldsymbol\sigma(s_i)t^{m_i-k}, \quad \text{if} \quad
\kappa_i=2m_i+1,
\; i=1,2, \, 
\end{aligned}
 \end{equation*}
  $$
  \sign (r):=
  \left \{
  \begin{array}{ll}
    0 & \, \text{if }\, r=0,\\
    1 & \, \text{if }\, r>0, \\
  \end{array}
\right .
  $$
and  $s_1=c$, $s_2=d$.
 \end{thm}

\textbf{Proof.} If $\kappa_c\geq 0$ and $\kappa_d\geq0$, the
operators $T(c)$ and $T(d)$ are right invertible. Moreover, the
corresponding operator $T(V(a,b))$ is also right-invertible and one
can easily check that the operator $U: H^p \times H^p \to H^p \times
H^p$,
  \begin{equation}\label{eqn17}
  U=
  \left(%
\begin{array}{cc}
 T_r^{-1}(c) T(\widetilde{a}^{-1})T_r^{-1}(d)  &   -T_r^{-1}(c)\\[1ex]
 T_r^{-1}(d)   & 0 \\
   \end{array}%
\right)
 \end{equation}
is one of the right-inverses for $T(V(a,b))$. Hence, one of the
solutions of equation \eqref{eqn13} is
 $$
(\Phi, \Psi)^T= U((2f,0)^T)=(2 T_r^{-1}(c)
T(\widetilde{a}^{-1})T_r^{-1}(d) f,2T_r^{-1}(d)f)^T,
 $$
and using formula \eqref{eqn9}, one obtains a solution  of the
non-homogeneous equation \eqref{eqn2}. Thus, if
$\kappa_c=\kappa_d=0$, the unique solution of the equation
\eqref{eqn2} is
 \begin{align*}
 \vp_0 &=  \frac 12 (\Phi-JQc\Phi+JQ\widetilde{a}^{-1}\Psi)\\
 &=T_r^{-1}(c) T(\widetilde{a}^{-1})T_r^{-1}(d) f -
JQcT_r^{-1}(c)T(\widetilde{a}^{-1})T_r^{-1}(d) f +JQ
\widetilde{a}^{-1} T_r^{-1}(d) f.
 \end{align*}
On the other hand, in the case where at least one of the indices
$\kappa_c$ or $\kappa_d$ does not vanish, the function $\vp_0$ is a
partial solution of \eqref{eqn2}. In order to obtain all solutions,
one has to add the solutions of the corresponding homogeneous
equation. However, the kernels of the operators $T(a)+H(b)$ have
been described earlier. Thus using Proposition 3.7 and Theorem 5.4
of \cite{DS:2014}, one arrives at the representation \eqref{eqn16}.
 \rbx

Now let us consider the situation where $\kappa_c$ and $\kappa_d$
are non-positive.
  \begin{thm}\label{t2}
Let $\kappa_c\leq 0$ and $\kappa_d\leq0$, and let $c_-$, $d_-$ be
factorization factors in the Wiener--Hopf factorization of the
matching functions $c$ and $d$, respectively.  If function $f\in
H^p$ satisfies the conditions
 \begin{equation}\label{eqn18}
  \begin{aligned}
\int_\sT \overline{d_-^{-1}(t)} t^j \overline{f(t)}\, |dt| & =0,
\quad j=0,1, \ldots, -\kappa_d-1 ,\\
\int_\sT  T_r^{-1}(\overline{d})
T(\overline{\widetilde{a}^{-1}})\left (\overline{c_-^{-1}(t)}
t^j\right ) \overline{f(t)}\, |dt| & =0, \quad j=0,1, \ldots,
-\kappa_c-1,
\end{aligned}
 \end{equation}
then the equation \eqref{eqn2} is uniquely solvable and  its
solution $\vp_0$ has the form
 \begin{multline}\label{eqn19}
\vp_0 =T_l^{-1}(c) T(\widetilde{a}^{-1})T_l^{-1}(d) f -
JQcT_l^{-1}(c)T(\widetilde{a}^{-1})T_l^{-1}(d) f +JQ
\widetilde{a}^{-1} T_l^{-1}(d) f .
\end{multline}
\end{thm}
\textbf{Proof.} If indices $\kappa_c$ and $\kappa_d$ are
non-positive, the operator $T(V(a,b))$ is left-invertible and
  \begin{equation*}
  T_l^{-1}(V(a,b))=
   \left(%
\begin{array}{cc}
 T_l^{-1}(c) T(\widetilde{a}^{-1})T_l^{-1}(d)  &   -T_l^{-1}(c)\\[1ex]
 T_l^{-1}(d)   & 0 \\
   \end{array}%
\right)
 \end{equation*}
is one of its left inverses. Therefore, equation \eqref{eqn13} is
solvable if and only if its right-hand side $(2f,0)^T$ is orthogonal
to any solution of the equation
 \begin{equation}\label{eqn20}
 T^*(V(a,b)) \psi =0,
 \end{equation}
where
 $$
T^*(V(a,b))=\left(%
\begin{array}{cc}
 0  & -T(\overline{c}) \\
  T(\overline{d})  & T(\overline{\widetilde{a}^{-1}}) \\
   \end{array}%
\right)
 $$
is the adjoint operator for the operator $T(V(a,b))$. However,
according to \cite[Proposition 3.3]{DS:2014}, the kernel of this
operator can be represented in the form
 $$
\ker
T^*(V(a,b))=\Omega(\overline{d})\dotplus\widehat{\Omega}(\overline{c}),
 $$
where
 \begin{align*}
\Omega(\overline{d}) &:=\left \{   (v,0)^T:v\in \ker T(\overline{d})\right\},\\
\widehat{\Omega}(\overline{c}) &:=\left \{ (T_r^{-1}(\overline{d})
T(\overline{\widetilde{a}^{-1}})u,u)^T:u\in \ker T(\overline{c})
\right \}.
\end{align*}
Taking into account the fact that
 \begin{align*}
    \ker T(\overline{d})& =\{\overline{d_-^{-1}}\, t^j: j=0,1,\ldots, -\kappa_d-1
    \}, \\
 \ker T(\overline{c})& =\{\overline{c_-^{-1}}\, t^j: j=0,1,\ldots, -\kappa_c-1
    \},
\end{align*}
one obtains solvability conditions \eqref{eqn12}. If equation
\eqref{eqn13} is solvable, its solution $\psi$ is
 $$
\psi=T_L^{-1}(V(a,b))(2f,0)^T=(2 T_l^{-1}(c)
T(\widetilde{a}^{-1})T_l^{-1}(d) f,2T_l^{-1}(d)f)^T,
 $$
and formula \eqref{eqn9} leads to the representation \eqref{eqn19}.
\rbx

  \begin{rem}\label{rem1}
Using Theorem 5.4 and Theorem 6.4 of \cite{DS:2014}, one can obtain
conditions of simultaneous solvability of the equations
\begin{align*}
    (T(a)+H(b))\vp&=f, \\
    (T(a)-H(b))\psi&=f,
\end{align*}
and show that if these two equations are solvable, then the
conditions \eqref{eqn18} are satisfied, so the system \eqref{eqn13}
is also solvable.
  \end{rem}

Consider one more case where all solution of the equation
\eqref{eqn2} can be found.

  \begin{thm}\label{t3}
Let $\kappa_c>0$, $\kappa_d<0$, and let $c_-$, $d_-$ be the
factorization factors in the Wiener--Hopf factorization of the
matching functions $c$ and $d$, respectively.  If function $f\in
H^p$ satisfies the conditions
 \begin{equation}\label{eqn21}
  \int_\sT \overline{d_-^{-1}(t)} t^j \overline{f(t)}\, |dt|  =0,
\quad j=0,1, \ldots, -\kappa_d-1 ,
 \end{equation}
then the equation \eqref{eqn2} is solvable  and all solutions of
this equation are given by the formula
  \begin{equation}\label{eqn22}
  \begin{aligned}
\vp &=T_r^{-1}(c) T(\widetilde{a}^{-1})T_l^{-1}(d) f -
JQcT_r^{-1}(c)T(\widetilde{a}^{-1})T_l^{-1}(d) f\\
&\quad +JQ \widetilde{a}^{-1} T_l^{-1}(d) f
 + c_+^{-1}\sum_{k=0}^{\kappa(1)} r_k^{(1)} u_k^{(\kappa(1),-)},
\end{aligned}
 \end{equation}
where $\kappa(1), r_k^{(1)}$, and $u_k^{(\kappa(1),-)}$ are defined
in Theorem \ref{t1}.
\end{thm}
\textbf{Proof.} Taking into account the conditions imposed on the
indices $\kappa_c$ and $\kappa_d$, one obtains that the operators
$T(c)$ and $T(d)$ are respectively, right- and left-invertible.
Further, straightforward computations show that the operator
   \begin{equation*}
  T_g^{-1}(V(a,b))=
   \left(%
\begin{array}{cc}
 T_r^{-1}(c) T(\widetilde{a}^{-1})T_l^{-1}(d)  &   -T_r^{-1}(c)\\[1ex]
 T_l^{-1}(d)   & 0 \\
   \end{array}%
\right),
 \end{equation*}
is a generalized inverse for the operator $T(V(a,b))$. Hence,
equation \eqref{eqn17} is solvable if and only if its right hand
side $(2f,0)^T$ is orthogonal to all solutions of the equation
\eqref{eqn20}. On the other hand, the operator $T^*(V(a,b))$ admits
the factorizations.
 \begin{align}
T^*(V(a,b))&
=\left(%
\begin{array}{cc}
 0  & -T(\overline{c}) \nn \\
  T(\overline{d})  & T(\overline{\widetilde{a}^{-1}}) \\
   \end{array}%
\right)\\[1ex]
 &= \left(%
\begin{array}{cc}
 T(\overline{c})  & 0 \\
 0  & I \\
   \end{array}%
\right)
\left(%
\begin{array}{cc}
 0  & -I \\
  I  & T(\overline{\widetilde{a}^{-1}}) \\
   \end{array}%
\right)
\left(%
\begin{array}{cc}
 -T(\overline{d})  &  0\\
  0  & I \\
   \end{array}%
\right).  \label{eqn23}
\end{align}
Note that the first operator in \eqref{eqn23} is invertible from the
left, whereas the operator
 \begin{equation*}
   D:=\left(%
\begin{array}{cc}
 0  & -I \\
  I  & T(\overline{\widetilde{a}^{-1}}) \\
   \end{array}%
\right)
\end{equation*}
is just invertible and
 \begin{equation*}
   D^{-1} =\left(%
\begin{array}{cc}
 T(\overline{\widetilde{a}^{-1}})  & I \\
 - I  & 0 \\
   \end{array}%
\right).
\end{equation*}
Therefore, factorization \eqref{eqn23} shows that the kernel of the
operator $T^*(V(a,b)$ can be represented in the form
 $$
\ker T^*(V(a,b)=\{(v,0)^T;v\in \ker T(\overline{d})\},
 $$
and the solvability condition \eqref{eqn21} follows. A partial
solution of the equation \eqref{eqn13} can be now obtained, viz.
$$
(\Phi, \Psi)^T= T_g^{-1}(V(a,b))(2f,0)^T,
$$
which allows one to construct a partial solution of the equation
\eqref{eqn2}. Adding the solutions of the corresponding homogeneous
equation, one obtains formula \eqref{eqn22}. \rbx

Consider now the last case remaining--viz. the situation where
$\kappa_c<0$ and $\kappa_d>0$. The method we use is not directly
applicable here. However, the initial equation admits a modification
which can be employed to find an analytic solution of the equation
\eqref{eqn2}. Before we go on, let us formulate a few auxiliary
results.

 \begin{lem}\label{l3}
 Let $f\in H^p$ and $n\in\sN$. If the equation
  \begin{equation}\label{Eq1}
  (T(a)+H(b))\vp =f
 \end{equation}
is solvable and $\vp_0$ is a solution of \eqref{Eq1}, then the
equation
  \begin{equation}\label{Eq2}
  (T(t^{-n}a) + H(t^n b))\psi =f
 \end{equation}
is also solvable and there is a solution $\psi_0$ of \eqref{Eq2}
which belongs to the image $\im T(t^n)$ of the operator
$T(t^n):H^p\to H^p$.
 \end{lem}
The proof is straightforward and is based on the relation
  \begin{equation}\label{Eq3}
  (T(a)+H(b))=(T(t^{-n}a) + H(t^n b))T(t^n).
 \end{equation}
Note that one of the solutions in question is $\psi_0=T(t^n)\vp_0$.
It is also worth mentioning that
 $$
\im T(t^n)=\{\psi\in H^p:
\widehat{\psi}_0=\widehat{\psi}_1=\cdots=\widehat{\psi}_{n-1}=0\},
 $$
where $\widehat{\psi}_j$, $j=0,1,\ldots, n-1$ are the Fourier
coefficients of the function $\psi$.
 \begin{cor}\label{c1}
Let $f\in H^p$ and $n\in\sN$. If equation \eqref{Eq2} does not have
any solutions $\psi\in \im T(t^n)$, then equation \eqref{Eq1} is not
solvable.
 \end{cor}

Thus any solution of \eqref{Eq1} produces a solution of \eqref{Eq2}
which lies in the set $\im T(t^n)$. On the other hand, one also has
a reverse statement.
  \begin{lem}\label{l4}
Let $f\in H^p$ and $n\in\sN$. If the equation \eqref{Eq2} is
solvable and has a solution $\psi_0\in \im T(t^n)$, then equation
\eqref{Eq1} is also solvable and one of its solutions has the form
 \begin{equation}\label{Eq5}
 \vp_0=T(t^{-n})\psi_0.
 \end{equation}
 \end{lem}
 \textbf{Proof.}
Use the fact that $T(t^n)T(t^{-n}):\im T(t^n)\to \im T(t^n)$ is the
identity operator. \rbx

Now we can describe how to find the solutions of equation
\eqref{Eq1} in the case at hand. Assume that this equation is
solvable and choose an $n\in \sN$ such that
\begin{equation*}
1\geq 2n+\kappa_c\geq 0.
\end{equation*}
Such an $n$ is uniquely defined and
\begin{equation*}
2n+\kappa_c =\left\{%
\begin{array}{ll}
   0, & \hbox{if\;} \kappa_c \; \hbox{is even,}  \\
    1, &\hbox{if\;} \kappa_c \; \hbox{is odd.} \\
\end{array}%
\right.
\end{equation*}
According to Lemma \ref{l3} equation \eqref{Eq2} is also solvable.
One can also see that $(at^{-n}, bt^n)$ is again a matching pair
with the subordinated pair $(ct^{-2n}, d)$. Moreover, the indices of
the corresponding operators $\widetilde{\kappa}_1=\ind T(c
t^{-2n})=\kappa_c+2n$ and $\widetilde{\kappa}_2=\ind T(d)=\kappa_d$
are non-negative. Thus the operator $T(at^{-n})+H(bt^n)$ satisfies
all conditions of Theorem \ref{t1}, and using the correspondingly
adapted formula \eqref{eqn16}, one obtains all solutions of the
equation \eqref{Eq2}. Since equation \eqref{Eq1} is solvable, the
set of solutions of \eqref{Eq2} contains at least one solution
$\psi_0$ which belongs to $\im T(t^n)$. Now one can employ formula
\eqref{Eq5} to obtain a solution of \eqref{Eq1}. The set of the
solutions of the corresponding homogeneous equation
 $$
(T(a)+H(b))\vp=0
 $$
is described in Theorem 6.3 of \cite{DS:2014}. This description
allows one to find all solutions of the equation \eqref{Eq1} in the
case $\kappa_c<0$, $\kappa_d>0$.

 \section{Examples}

Let us illustrate the above theory by a few simple examples.

 \begin{exm}\label{ex1}
Consider the equation
 \begin{equation}\label{eqn24}
( T(t^{-2})+H(t^2))\vp(t)=f(t).
 \end{equation}
Obviously, the functions $a(t)=t^{-2}$ and $b(t)=t^2$ constitute a
matching pair. Moreover, one has
\begin{equation*}
 \begin{array}{lll}
\widetilde{a}^{-1}(t)=t^{-2}, & \quad c(t)=a(t) b^{-1}(t)=t^{-4},
&\quad d(t)=
\widetilde{a}^{-1}(t) b(t)=1,\\[1.5ex]
   c_+(t)=c_-(t)=1, & \quad
 d_+(t)=d_-(t)=1,& \quad
 \boldsymbol\sigma(c)=\boldsymbol\sigma(d)=1.\\[1.5ex]
  \kappa_c=4, &\quad \kappa_d =0. &
\end{array}
 \end{equation*}
Thus equation \eqref{eqn24} is subject to Theorem \ref{t1}. It is
solvable for any right-hand side $f\in H^p$, $1<p<\infty$ and any
solution of \eqref{eqn24} can be obtained from formula
\eqref{eqn16}. Taking into account that $T_r^{-1}(c)=Pt^4P$,
$T^{-1}(d)=I$, the solutions of this equation can be written in the
form
 \begin{equation}\label{eqn25}
 \vp (t) = Pt^4Pt^{-2}Pf(t) - JQt^{-4}Pt^4Pt^{-2}Pf(t) +JQ t^{-2}
 f(t) + r_1(t-t^2) +  r_2(t^2-t^3),
 \end{equation}
where $r_1,r_2$ are arbitrary complex numbers. In particular, let us
find a solution of the equation \eqref{eqn24} for a given right-hand
side $f$. For example, if $f(t)=t^6+3t^4$ and $r_1=0$, $r_2=0$, the
formula \eqref{eqn25} produces the function
 $$
\vp(t)=t^8+3 t^6,
 $$
and one can easily check that this is a partial solution of the
equation \eqref{eqn24} with the right-hand side $f(t)=t^6+3t^4$.
 \end{exm}

  \begin{exm}\label{ex2}
Consider the equation
 \begin{equation}\label{eqn26}
( T(2t+1)+H(2t+1))\vp(t)=f(t).
 \end{equation}
It is clear that $(2t+1,2t+1)$ is a matching pair. Further, one has
\begin{equation*}
\widetilde{a}^{-1}(t)=\frac{t}{t+2}, \quad  c(t)=1, \quad
d(t)=\frac{t(2t+1)}{t+2},
 \end{equation*}
and the function $d(t)$ admits the factorization
 $$
d(t)= \left (\frac{2}{t+2} \right )\, t^2 \,\left
(\frac{2t+1}{2t}\right )\,,
 $$
with the factorization factors
 $$
d_+(t)=\frac{2}{t+2}, \quad d_-(t)=\frac{2t+1}{2t}.
 $$
Thus
 $$
 \kappa_c=0, \quad \kappa_d =-2,
 $$
and the equation \eqref{eqn26} is subject to Theorem \ref{t2}.
Therefore, on order to use the above method, we need the right-hand
side $f$ to satisfy the condition
 \begin{equation}\label{eqn27}
 \int_\sT \frac{t^j}{t+2} \overline{f(t)}\, |dt|  =0,
\quad j=0,1.
 \end{equation}
If these conditions are satisfied, then equation \eqref{eqn26} is
uniquely solvable and its solution can be found by the formula
\eqref{eqn19}. Consider, for example, the function
$f(t)=(2t+1)(t^2+t)$. One can easily check that this function $f$
satisfies the solvability conditions \eqref{eqn27}. Applying formula
\eqref{eqn19} with
 $$
T_l^{-1}(d)= Pt^{-2}P(t+2) P \left (\frac{t}{2t+1} \right )P,
 $$
one obtains the function
 $$
\vp_0 (t)=t(t+1),
 $$
which is the solution of \eqref{eqn26} with the right-hand side
$f(t)=(2t+1)(t^2+t)$.

On the other hand, it is worth mentioning that in this case the
operator $T(V(a,b))$ is not right invertible. Therefore, there are
right-hand sides $f\in H^p$ such that the equation \eqref{eqn26} is
solvable but its solution cannot be found by the method used.
 \end{exm}

\begin{exm}\label{ex3}
Let $b=b(t)$ be a matching function, that is $b(t)
\widetilde{b}(t)=1$. Consider the following equation
  \begin{equation}\label{Eq6}
  \vp + H(b)\vp=f, \quad f \in H^p.
 \end{equation}
Assume that the operator $T(b)$ is Fredholm. Then according to
\eqref{Eq27} the function $b$ admit  the Wiener--Hopf factorization
of the form
\begin{equation}\label{Eq7}
b(t)=b_+(t) t^n (\boldsymbol\sigma(b) \widetilde{b}_+^{-1}(t)).
 \end{equation}
One also has $a(t)=1, t\in \sT$, so that
 $$
 c(t)=b^{-1}(t)=\widetilde{b}(t), \quad d(t)=b(t).
 $$
Therefore, $d$ has the factorization \eqref{Eq7} and $c$ can be
factorized as follows
 $$
c(t)= b_+^{-1}(t) t^{-n}(\boldsymbol\sigma(b) \widetilde{b}_+(t)).
 $$
 Thus
  $$
\kappa_c=n, \quad \kappa_d=-n,
  $$
and if $n\neq0$,  the indices $\kappa_c$ and $\kappa_d$ have
different signs. Assume for definiteness that $n\geq 0$. Then we are
in the situation described by Theorem \ref{t3}, and if
 $$
\int_\sT \overline{\widetilde{b}}_+ (t) \overline{f(t)} \,t^j \,
|dt| = 0, \quad j=0,1,\ldots, n-1,
 $$
one can write the solutions in the form \eqref{eqn22}. In
particular, one has
 \begin{align*}
T_r^{-1}(c)& =\boldsymbol\sigma(b) T(b_+)
T(\widetilde{b}_+^{-1})T(t^n)\\
T_l^{-1}(d)& =\boldsymbol\sigma(b)T(t^{-n}) T(b_+^{-1})
T(\widetilde{b}_+).
\end{align*}
 Moreover, let $Q_n$ be the operator defined by
  $$
Q_n \vp(t)= Q_n\left( \sum_{j=0}^\infty \widehat{\vp}_j t^j\right
):=\sum_{j=n+1}^\infty \widehat{\vp}_j t^j.
  $$
One can easily see that $Q_n=T(t^n)T(t^{-n})$. Therefore, all
solutions of the equation \eqref{Eq6} can be written in the form
 $$
\vp= (I-JQb^{-1})  T(b_+) T(\widetilde{b}_+^{-1}) Q_n T(b_+^{-1})
T(\widetilde{b}_+) f +b_+\sum_{j=0}^{\kappa(1)} r_j u_j^{(k(1),-)},
 $$
 where $r_j, j=0,1,\ldots, \kappa(1)$ are arbitrary complex numbers.
  \end{exm}


\end{document}